\documentstyle[12pt]{article}
\newcommand{\const}{\mathop{\rm const}\limits}

\newcommand{\grad}{\mathop{\rm grad}\limits}

\newcommand{\supp}{\mathop{\rm supp}\limits}

\newcommand{\Rad}{\mathop{\rm Rad}\limits}

\begin{document}
\begin{center}

{\bf Exact constant in Sobolev's  and Sobolev's trace inequalities }\\

\vspace{3mm}

  {\bf for Grand Lebesgue Spaces with monomial weight. }\\

\vspace{3mm}

{\sc Ostrovsky E., Sirota L.}\\

\normalsize

\vspace{3mm}
{\it Department of Mathematics and Statistics, Bar-Ilan University,
59200, Ramat Gan, Israel.}\\
e-mail: \ eugostrovsky@list.ru \\

\vspace{3mm}

{\it Department of Mathematics and Statistics, Bar-Ilan University,
59200, Ramat Gan, Israel.}\\
e-mail: \ sirota3@bezeqint.net \\

\end{center}

\vspace{4mm}

 {\it Abstract.}  We generalize in this article the classical Sobolev's  and Sobolev's trace inequalities  on the Grand
 Lebesgue Spaces under monomial weight instead the classical Lebesgue or grand Lebesgue Spaces. \par
 We will distinguish the classical Sobolev's inequality and the so-called trace Sobolev's inequality. \par

 \vspace{3mm}

 {\it Key words:} Sobolev's and Sobolev's trace inequalities, derivative, gradient, norm, monomial and monomial weight, dilation operator,
 Lebesgue spaces, Talenti's estimate, Bilateral  Grand Lebesgue  spaces, scaling, fundamental function. \par

 \vspace{3mm}

{\it Mathematics Subject Classification (2000):} primary 60G17; \ secondary
 60E07; 60G70.\\

\vspace{3mm}

\section{Introduction. Notations. Statement of problem.}

\vspace{3mm}

{\bf A. Classical Sobolev's inequality.} \par

 The classical Sobolev's inequality in the whole space $ R^m, $ see, e.g.
 \cite{Kantorovicz1}, chapter 11, section 5; \cite{Bliss1}, \cite{Sobolev1}, \cite{Talenti1} etc.
 asserts that for all function $ f, f: R^m \to R, \ m \ge 3 $ from the Sobolev's space
 $ W^1_p(R^m), $ which may be defined as a closure in the Sobolev's norm

 $$
 ||f||W^1_p(R^m)= |f|_p + |Df|_p
 $$
  of the set of all finite continuous differentiable functions $ f, f: R^m \to R, $ that

  $$
  |f|_{q} \le K_m(p) \ |Df|_p, \ q = q(p)= mp/(m-p), \ p \in [1,m), \
  q \in (m/(m-1), \infty).\eqno(1.1)
  $$
  Here  $ m = 3,4, \ldots; $

  $$
  |f|_p = |f|_{p,m} = |f|_{p,R^m} = \left[ \int_{R^m} |f(x)|^p \ dx  \right]^{1/p},
  $$

  $$
  Df = \{ \partial f/\partial x_1, \partial f/\partial x_2, \partial f/\partial x_3,
\ldots,  \partial f/\partial x_m \} = \grad_x f,
  $$

$$
|Df|_p = \left| \left[ \sum_{i=1}^m ( \partial f/\partial x_i)^2  )  \right]^{1/2} \right|_p.
$$
The best possible constant in the inequality (1) belongs to G.Talenti \cite{Talenti1}, see also  \cite{Bliss1}:

$$
K_m(p) = \pi^{-1/2} m^{-1/p} \left[\frac{p-1}{m-p}\right]^{1-1/p} \cdot
\left[\frac{\Gamma(1+m/2) \ \Gamma(m)}{\Gamma(m/p) \ \Gamma(1+m-m/p)} \right]^{1/m}.
$$

\vspace{3mm}

{\bf B. Trace Sobolev's inequality.} \par

\vspace{3mm}

Let $ m,n = 1,2,\ldots, \ x \in R^m, y \in R^n, z = \{x,y \} \in R^{m+n}, \ u =
u(x,y) = u(z)  $ be any function from the space $ W^1_p(R^{m+n}). $ \par
 We consider in this case only the so-called {\it radial } functions. In detail, we
 define as usually for the vectors $ x = \vec{x} = \{x_1,x_2, \ldots, x_m \} $ and
 $ y = \vec{y} = \{y_1,y_2, \ldots, y_n \} $

 $$
 |x| = \left(\sum_{i=1}^m (x_i)^2 \right)^{1/2}, \ |y| =
 \left(\sum_{j=1}^n (y_j)^2   \right)^{1/2}
 $$
 and correspondingly

 $$
 |z| = |(x,y)| = \left(|x|^2 + |y|^2 \right)^{1/2}.
 $$



Let us denote $ N = m + n, (N \ge 3); \ S[u](x) = u(x,0),$
$$
 \nabla u = \{ \partial u/\partial x_1, \partial u/\partial x_2, \ldots,
 \partial u/\partial x_m, \partial u/\partial y_1, \partial u/\partial  y_2,
  \ldots,  \partial u/\partial y_n \}=
$$

$$
\{\grad_x u, \ \grad_y u \}; \ |\nabla u|_p = (|\grad_x u|^p_p +
 |\grad_y u|^p_p)^{1/p}.
$$

We will denote the class of all the radial functions $ \Rad = \Rad(R^N); \ u(\cdot)
\in \Rad. $ \par

 Notice that the operator $ S[u] $ is correct and continuously defined in the
 $ L_p(R^m) $ in the following sense:

 $$
\lim_{|y| \to 0} |u(\cdot,y) - S[u]|L_p(R^m) = 0,
 $$
see \cite{Besov1},  chapter 5, section 24.\par

The following inequality is called the {\it  Sobolev's trace inequality:}

$$
|S[u](\cdot)|_{q,m} \le K_{m,n}(p) \cdot |\nabla u|_{p,N}, \ q = q(p) = mp/(N-p), \
p \in [1,N). \eqno(1.2).
$$
   We will understand further under the constant $ K_{m,n}(p) $ in the inequality
(1.2) its minimal value, namely:

$$
K_{m,n}(p)= \sup
\left\{ \left[ \frac{|S[u](\cdot)|_q}{|\nabla u|_p} \right], \
u \in W^1_p(R^{m+n})\cap \Rad(R^N), \nabla u \ne 0 \right\}. \eqno(1.3)
$$

 It is evident $ K(m,0) = K(m). $ \par
More information about the constant $ K_{m,n}(p)$ see, for instance, in the articles
\cite{Beesack1}, \cite{Nazaret1}, \cite{Young1}, \cite{Biezuner1}, \cite{Druet1}, \cite{Lieb1}, \cite{Zhu1},
\cite{Edmunds1}, \cite{Escobar1}, \cite{Maz'ja1} etc., see also reference therein.\par

\vspace{4mm}

{\bf  We intend in this article to extend the Sobolev's inequalities into two directions: on the case of
monomial weight inequalities and into the so-called Weight Grand Lebesgue spaces (WGLS). } \par

\vspace{3mm}

{\bf We intend to show also the exactness of offered estimations.}

 \vspace{4mm}

 The {\it monomial} Sobolev's inequality (imbedding theorems)
  in the classical Lebesgue-Riesz spaces $ L_p(R^N) $ with exact constant
computation  is obtained in  the article  of Xavier Cabre and Xavier Ros-Oton
\cite{Cabre1}; see also reference therein. \par
Some preliminary  results see in \cite{Ostrovsky500} - \cite{Ostrovsky503}.\par
 Many interest applications of these inequalities in the theory of PDE and in the theory of function
are also described in the source article \cite{Cabre1}.\par

\vspace{3mm}

{\bf Remark 1.1.}  "However, in general the monomial
weight does not satisfy the Muckenhoupt condition $ A(p) $ and Theorem 1.3 cannot
be deduced from these results on weighted Sobolev inequalities, even without the
best constant in the inequality;" see \cite{Cabre1}. \par

\vspace{3mm}

 We  must introduce some new notations, following in the the majority
 \cite{Cabre1}. Let $ A = (A(1), A(2), \ldots, A(m)) $ be
 $  m \ - $ tuple of non-negative numbers, briefly: $ A \in T(m). $ The {\it monomial} $ x^A, \ x \in R^m, $
by definition, is a function of a form

$$
x^A = |x|^A =  |x_1|^{A(1)} \  |x_2|^{A(2)} \  |x_m|^{A(m)}  = \prod_{i=1}^m  |x_i|^{A(i)}. \eqno(1.4)
$$

 The correspondent monomial measure $  \mu_A  $ on the whole space $  R^m $ may be defined as follows:

$$
\mu_A(G) = \int_G x^A \ dx; \eqno(1.5)
$$

$$
\int f(x) dx \stackrel{def}{=} \int_{R^m} f(x) dx.
$$

\vspace{4mm}

{\bf Definition 1.1.}  Let $  A,B \in T(m) $ be some $ m \ $ tuples, $ p,q = \const > 0.$ If  there exists a finite positive
constant, more exactly, function $ K = K_{A,B}(p,q)  $ on $ (A,B; p,q) $ such that the inequality of the Sobolev's form

$$
\left( \int |x|^B \ |u(x)|^q \ dx   \right)^{1/q} \le K_{A,B}(p,q)  \cdot \left( \int |x|^A \ |\nabla u(x)|^p \ dx \right)^{1/p} \eqno(1.6)
$$
or equally $ L(u) = L_{A,B;p,q}(u) := $

$$
 |u(\cdot)|_{q, \mu_B}  \le K_{A,B}(p,q) \cdot |\nabla u(\cdot)|_{p, \mu_A} =: R(u) = R_{A,B;p,q}(u) \eqno(1.6a)
$$
there holds  for arbitrary infinite differentiable function $  u = u(x), \ x \in R^m $ with compact support: $  u(\cdot) \in C_0^{\infty}(R^m),  $
 then we will write $ (A,B; p,q) \in S $ or equally will talk that the four $ (A,B; p,q)  $ is Sobolev's four. \par

\vspace{4mm}

 Hereafter $ C, C_j $ will denote any non-essential finite positive constants.
 We define also for the values $ (p_1, p_2), $ where $ 1 \le p_1 < p_2 \le \infty $

$$
L(p_1, p_2) =  \cap_{p \in (p_1, p_2)} \ L_p.
$$

\vspace{3mm}



\vspace{4mm}

\section{ Necessary condition for Sobolev inequality with monomial weight}

\vspace{4mm}

 Denote

 $$
 D(A) = \sum_{i=1}^m A(i) + m, \ D(B) = \sum_{i=1}^m B(i) + m. \eqno(2.0)
 $$

\vspace{3mm}

{\bf Theorem 2.1.} {\it Suppose $ m \ge 1, \  D(A) > 1, \ D(B) > 0. $ If the four $  (A,B;p,q) $ is the Sobolev's four, then }

$$
q = \frac{D(B) \ p}{ D(A) - p}. \eqno(2.1)
$$
{\it As a  consequence: in this case} $ p < D(A). $ \par

\vspace{3mm}

{\bf Proof.} The relation (2.1) may be obtained by means of the so-called {\it dilation method}, or equally {\it scaling method},
 see \cite{Stein1},   \cite{Talenti1}. This method was used,  e.g., in
 \cite{Ostrovsky3},  \cite{Ostrovsky4}, \cite{Ostrovsky6}, \cite{Ostrovsky503}. \par

In detail, let us define as usually the family of dilation operators
$ T_{\lambda}[u], \ \lambda \in (0,\infty), \ u: R^m \to R, $ of a form:

$$
T_{\lambda}[u](x) = u(\lambda \cdot x). \eqno(2.2)
$$
 Obviously, if $ u(\cdot) \in C_0^{\infty}(R^m), $ then  $ T_{\lambda} u(\cdot) \in C_0^{\infty}(R^m). $ \par
Let $ u = u(x) $ be some non-zero function  from the set $ \in C_0^{\infty}(R^m) $ for which the inequality (1.6) is
satisfied. We deduce  substituting into (1.6a) the function $ u_{\lambda} = T_{\lambda}u $

$$
L_{\lambda}:= |u_{\lambda}(\cdot)|_{q, \mu_B}  \le K_{A,B}(p,q) \cdot |\nabla u_{\lambda}(\cdot)|_{p, \mu_A} =: R(u_{\lambda}) = R_{A,B;p,q}(u_{\lambda})
=: R_{\lambda}. \eqno(2.3)
$$
 It is easy to conclude

$$
L_{\lambda} = \lambda^{ -(m + \sum B(i))/q  } \ L =  \lambda^{ -D(B)/q  } \ L,
$$

$$
R_{\lambda} = \lambda^{ 1-(m + \sum A(i))/p  } \ R =  \lambda^{ 1-D(A)/p  } \ R.
$$
We have  substituting into (1.6a)

$$
\lambda^{ -D(B)/q  } \ L \le K_{A,B}(p,q)\ \lambda^{ 1-D(A)/p  } \ R. \eqno(2.4)
$$
 The last inequality (2.4) may be true for arbitrary positive value $ \lambda $ if and only if

$$
- D(B)/q = 1 - D(A)/p,
$$
which is  equivalent to the equality (2.1).\par

\vspace{3mm}

{\bf Remark 2.1.} In the case when $ A = B $ or at least if $ \sum B(i) = \sum A(i), $ the relation (2.1) coincides with the
{\it  sufficient} condition number "(1.5)" for the monomial Sobolev's inequality in the article \cite{Cabre1}. \par

\vspace{4mm}

\section{ Necessary condition for traces Sobolev's inequalities with monomial weight}

\vspace{4mm}

 Let us consider in this section the {\it trace} inequality Sobolev's type with monomial weight  of the form

$$
\left(\int_{R^r} |x|^B \ |u(x,0)|^q \ dx  \right)^{1/q} \le K^{Tr}_{A,B}(p,q) \cdot
\left( \int_{R^d} |\{x,y\}|^A \ |\nabla u(x,y)|^p \ dx dy  \right)^{1/p}=
$$

$$
K^{Tr}_{A,B}(p,q) \cdot \left( \int_{R^d} |z|^A \ |\nabla u(z)|^p \ dz  \right)^{1/p}. \eqno(3.1)
$$
 Here $  z = \vec{z} = \{ x,y \},  $
 $$
 \dim x = r, 1 \le r < d, \ d = \dim \{x,y\} = r + \dim y, \ B = (B(1), B(2), \ldots, B(r)),
 $$

 $$
  B(i) \ge 0, \  A = (A(1), A(2), \ldots, A(r), A(r+1), \ldots, A(d)), \ A(j) \ge 0;
 $$

$$
|x|^B = \prod_{i=1}^r |x_i|^{B(i)}, \ |\{x,y\}|^A = \prod_{j=1}^r |x_j|^{A(j)} \cdot \prod_{k=r+1}^d |y_k|^{A(k)}, \eqno(3.2)
$$
and we define

$$
D_r(B) = r + \sum_{i=1}^rB(i),
$$
so that $  D_m(A) = D(A). $\par

  We find as before:\par

  \vspace{4mm}

 {\bf Theorem 3.1.} {\it Assume  $ D(A) > 1. $
 If the inequality (3.1) holds true for arbitrary function $ u(\cdot) $  from the space $ C^{\infty}_0(R^d),  $ then}

$$
q = \frac{D_r(B) \ p}{D(A) - p}. \eqno(3.3)
$$
{\it  As a consequence: in the considered case as before } $  p < D(A). $ \par

\vspace{4mm}

\section{Grand Lebesgue Spaces.}\par

 {\bf Definition.} \par

\vspace{3mm}

  Recently, see \cite{Kozachenko1}, \cite{Fiorenza1}, \cite{Fiorenza2}, \cite{Fiorenza3}, \cite{Iwaniec1}, \cite{Iwaniec2},
 \cite{Ostrovsky1}, \cite{Ostrovsky2}, \cite{Ostrovsky3}, \cite{Ostrovsky4},  \cite{Ostrovsky5}, \cite{Ostrovsky6}  etc.
 appears the so-called {\it Grand Lebesgue Spaces} $ GLS = G(\psi) =  G(\psi; a,b), \ a,b = \const, a \ge 1, a < b \le \infty, $ spaces consisting
 on all the measurable functions $ f: T \to R, $ where $ (T =  \{t \}, M, \nu) $ is measurable space with non-trivial sigma-finite measure $  \nu, $
 having finite norms

$$
  ||f||G(\psi) = ||f||G(\psi; \nu)  \stackrel{def}{=} \sup_{p \in (a,b)} \left[ |f|_{p,\nu} /\psi(p) \right], \eqno(4.1)
$$

$$
|f|_{p,\nu} := \left[ \int_T |f(t)|^p \ \nu(dt)  \right]^{1/p}.
$$
 Here $ \psi(\cdot) $ is some continuous positive on the {\it open} interval  $ (a,b) $ function such that

$$
     \inf_{p \in (a,b)} \psi(p) > 0.
$$

We can suppose without loss of generality

$$
\inf_{p \in (a,b)} \psi(p) = 1.
$$

 Notation: $  (a,b) = \supp \psi. $\par

 As the capacity of the measure $ \nu $ may be picked the measures $ \mu_A $ or $  \mu_B $ in the case of course when
$  T = R^m $  equipped with Borelian sigma field $  M. $ In this case we will denote for brevity

$$
||f||G(\psi; \mu_A) = ||f||_A(G(\psi)). \eqno(4.2)
$$

    This spaces are rearrangement invariant, see \cite{Bennet1}, and
  are used, for example, in the theory of probability \cite{Talagrand1}, \cite{Kozachenko1}, \cite{Ostrovsky1};
  theory of Partial Differential Equations \cite{Fiorenza2}, \cite{Iwaniec2}; functional analysis \cite{Ostrovsky4},
  \cite{Ostrovsky5}; theory of Fourier series \cite{Ostrovsky7}, theory of martingales \cite{Ostrovsky2} etc.\par

 Let $ \delta = \const > 0; $ the fundamental function  $ \phi(\delta) = \phi_{ G(\psi)}(\delta)  $  of the space $  G(\psi) $
may be calculated as follows:

$$
\phi_{ G(\psi)}(\delta) = \sup_{p \in (a,b)} \left[ \frac{\delta^{1/p}}{\psi(p)} \right]. \eqno(4.3)
$$
 This function play a very important role in the theory of Fourier series, operator theory etc., see \cite{Bennet1}.
We intend to use further this notion in the theory of Sobolev's monomial inequalities. \par

\vspace{3mm}

 Suppose now the measure $ \nu $ in (4.1) is probabilistic: $ \nu(T) = 1, $
and let the function $ \psi = \psi(p) $  be such that $ \supp \psi = (1,\infty). $ Then the  Grand Lebesgue Space $ G(\psi) $
coincides up to norm equivalence to the subspace of all mean zero:  $ \int_T f(t) \ \nu(dt) = 0  $
 measurable function (random variables) of the so-called {\it exponential Orlicz space  } $ L(N) = L(N; T,\nu) $
 with {\it  exponential } Orlicz-Young function $ N = N(u),$ and conversely proposition is also true: arbitrary exponential
Orlicz space $ L(N) $  coincides with some Grand Lebesgue Space, see \cite{Kozachenko1}.\par
 We intend to obtain in the next section the Sobolev's inequalities with monomial weight for Grand Lebesgue Spaces; thus,
we will obtain as a slight  generalization these inequalities for exponential Orlicz spaces due in the article
\cite{Cabre1}.\par

\vspace{3mm}

  \section{Sobolev's inequality for weighted Grand Lebesgue Spaces.}

\vspace{3mm}

 We assume during this section in the monomial Sobolev's inequality $  A  = B, \ 1 \le p < D(A). $\par

Some new notations. As long as in this section $ A = B, $ we denote for simplicity $ D = D(A); $ recall that we suppose
 $  D > 1. $\par

 Further, we denote following \cite{Cabre1}

$$
C_1 = D \cdot \left( \frac{\Gamma( (1+ A(1))/2 )\ \Gamma( (1+ A(2))/2 ) \ \ldots \Gamma( (1+ A(m))/2 )}{2^k \ \Gamma( (1 + D)/2) } \right)^{1/D},
\eqno(5.1)
$$

$$
C(p) = C_1 \cdot D^{1/D - 1 - 1/p} \cdot \left( \frac{p-1}{D-p}  \right)^{1/p'} \cdot
\left( \frac{p' \ \Gamma(D)}{\Gamma(D/p) \ \Gamma(D/p')}    \right)^{1/D}, \eqno(5.2)
$$
where $ p' = p/(p-1), \ 1 < p < D  $ and  $ k $ is the number of positive  entries in the vector $ A. $ \par

 It is easily to verify that $ C(p)  $ is continuous function in a {\it semi - closed segment} $ 1 \le p < D:  $

$$
\lim_{p \to 1+0} C(p) = C_1,
$$
and  that as $ p \to D-0  $

$$
C(p) \sim C_1 \ (D - p)^{ -(1 - 1/D)}.
$$

 This constant (more exactly, function on $  p ) $ are the {\it exact} constants in the monomial Sobolev's inequality

$$
\left( \int |x|^A \ |u(x)|^q \ dx   \right)^{1/q} \le C(p)  \cdot \left( \int |x|^A \ |\nabla u(x)|^p \ dx \right)^{1/p}, \eqno(5.3)
$$
where the expected value $ q  $ is following:

 $$
 q = D \ p /(D - p); \eqno(5.4)
 $$
obtained in the article  \cite{Cabre1}  by means of isoperimetric  inequalities .\par

 "It is a surprising fact that the weight $ x^A $ is not radially symmetric but still Euclidean
balls centered at the origin (intersected with $ R^n) $ minimize this isoperimetric quotient", \cite{Cabre1}. \par

 The equality (5.4) defines uniquely the value $ p $ as a function on $  q: \ p = p(q) = qD/(q + D) $ and inversely $  q  $ as a function on $ p. $
Herewith

$$
1 \le p  < D \ \Leftrightarrow \frac{D}{D-1} < q  < \infty.
$$

\vspace{4mm}

 Let $ \psi(\cdot) $ be arbitrary function from the set $G\Psi(1,D). $ We define the new
function  $ \zeta_{\psi}(q) = \zeta(q) $ from the set $ G\Psi(D/(D-1), \infty) $ as follows:

$$
\zeta_{\psi}(q) = C \left( \frac{Dq}{D+q} \right) \cdot
\psi\left( \frac{Dq}{D+q} \right), \ q \in [D/(D-1), \infty). \eqno(5.5)
$$

{\bf Theorem 5.1}. {\it The following Sobolev's type inequality holds:}

$$
||u||G(\zeta_A) \le 1 \cdot ||\nabla u||G(\psi_A), \eqno(5.6)
$$
{\it  and the constant "one" in the inequality (5.6) is the best possible.} \par

\vspace{4mm}

{\bf Remark 5.1.} It is presumed in the assertion (5.6) of theorem 5.1 that the right-hand
 side of the proposition of theorem 5.1 is finite. \par

 \vspace{4mm}

{\bf Proof.} \par
{\bf A. The upper bound.} Let the function $ u = u(x) $ be such that

$$
\nabla u \in \cap_{p \in (1,D)} W^1_p(R^{m})
$$
and
$$
||\nabla u||G(\psi_A) \in (0,\infty);
$$
we can assume without loss of generality
$$
||\nabla u||G(\psi_A) = 1.
$$
It follows from the direct definition of norm for the $ G(\psi) $ spaces that

$$
|\nabla u|_{p, \mu(A)} \le \psi(p), \ p \in (1,D).
$$

We use  now the inequality (5.3):

$$
|u|_{q,A} \le C(p) \cdot \psi(p), \ p \in (1,D). \eqno(5.7)
$$
 Since the inequality (5.7) is valid for all the values inside the whole segment $ (1,D), $ it may be rewritten
taking into account the monotonicity of the function $ q  \to p(q) $  as follows:

 $$
|u|_{q,A} \le C(p(q)) \cdot \psi(p(q)) = \zeta(q), \ q \in (D', \infty). \eqno(5.8)
$$

The proposition of theorem 5.1 follows immediately from (5.8) after the substitution
$$
p = \frac{Dq}{D+q}, \ q \in (D(D-1), \infty)
$$
by virtue of the direct definition of the norm in GLS:

$$
||u||G(\zeta_A) \le 1  = 1  \cdot ||\nabla u||G(\psi_A). \eqno(5.9)
$$

\vspace{3mm}

{\bf B. Proof of the low bound.}\par

 The case $  A = B = 0 $ is provided in \cite{Ostrovsky501} by means of construction of corresponding example;
the general case follows immediately from the main result of the  article  \cite{Ostrovsky502}. \par

\vspace{3mm}

\section{Trace Sobolev's inequality for monomial Grand Lebesgue Spaces. \ Radial case.}

\vspace{3mm}

We consider in this section {\it only radial } functions $ u(z) = u(x,y) = u(|z|). $  Here

$$
q = \frac{p \ D_r(B)}{ D(A) - p}  \stackrel{def}{=} \frac{p \ D_r}{ D - p},  \ 1 \le p < D;
$$
and we consider in this section the {\it trace} inequality Sobolev's type with monomial weight  of the form

$$
\left(\int_{R^r} |x|^B \ |u(x,0)|^q \ dx  \right)^{1/q} \le K^{Tr}_{A,B}(p,q) \cdot
\left( \int_{R^d} |\{x,y\}|^A \ |\nabla u(x,y)|^p \ dx dy  \right)^{1/p}=
$$

$$
K^{R}_{A,B}(p q), \cdot \left( \int_{R^d} |z|^A \ |\nabla u(z)|^p \ dz  \right)^{1/p}, \
K^{R}_{A,B}(p q) \stackrel{def}{=} K^{R}. \eqno(6.1)
$$
for radial functions. \par

 Recall that the Jacobian of the spherical  $ r \ -  $ dimensional
coordinates relative the ordinary Euclidean ones $ J_r(\rho,\theta), \theta \in \Theta_r  $
has a form

$$
 J_r(\rho,\theta) = \rho^{r-1} \Omega_r(\theta),
$$
 and analogously

$$
 J_d(\rho,\phi) = \rho^{d-1} \Omega_d(\phi), \ \phi \in \Theta_d.
$$
 Define the following integrals:

$$
\omega_r^q(B) := \int_{\Theta_r} \Omega_r^{B + 1}(\theta) \ d \theta, \
\omega_d^p(A) = \int_{\Theta_d} \Omega_r^{A + 1}(\phi) \ d \phi, \eqno(6.2)
$$

$$
A + 1 = \{A(1) + 1, \ A(2) + 1, \ \ldots, A(d) + 1 \}, \  B + 1 = \{B(1) + 1, \ B(2) + 1, \ \ldots, B(r) + 1 \}.
$$
 The integrals (6.2)  are in fact calculated in \cite{Cabre1}. \par

 For the radial  functions $  u = u(|z|) $ the inequality (6.1) nay be transformed  as follows:

$$
\left[ \int_0^{\infty} \rho^{D_r -1} \ |u(\rho)|^q \ d \rho \right]^{1/q} \le \omega_d(A) \ \omega_r^{-1} (B) \ K^R_{A,B}(p,q) \times
$$

$$
\left[ \int_0^{\infty} \rho^{ D - 1} \ |u'(\rho)|^p  \ d \rho  \right]^{1/p}, \eqno(6.3)
$$
or equally

$$
\left[\int_0^{\infty} s^{D_r - 1} ds \cdot \left|  \int_s^{\infty} g(t) \ dt  \right|  \right]^{1/q} \le W_{A,B}(p,q; r,d) \cdot
\left[ \int_0^{\infty} |g'(s)|^p \ ds  \right]^{1/p}, \eqno(6.3a)
$$
where

$$
W_{A,B}(p,q; r,d) = \omega_d(A) \ \omega_r^{-1} (B) \ K^R_{A,B}(p,q).
$$
 Moreover, we will understand under the value $ W_{A,B}(p,q; r,d) $ its optimal, i.e. minimal value. \par
 Denote also

$$
M = M_{A,B}(p,q; r,d) = D_r^{-1/r} \cdot \left( \frac{p-1}{D-r} \right)^{1 - 1/p}, \
$$

$$
Q = Q_{A,B}(p,q; r,d) = \left( \frac{ q}{q-1}  \right)^{1  - 1/p}  \cdot q^{1/q}.
$$

 \vspace{3mm}

{\bf Theorem 6.1.} {\it Under formulated above conditions }

$$
M_{A,B}(p,q; r,d) \le W_{A,B}(p,q; r,d) \le M_{A,B}(p,q; r,d) \cdot Q_{A,B}(p,q; r,d).    \eqno(6.4)
$$

\vspace{3mm}

{\bf Proof.}
 It is sufficient to use the result belonging to Bradley \cite{Bradley1};
 see also the famous monograph of Maz'ja \cite{Maz'ja1}, p. 45, Theorem 3. This estimate is a weight generalization of the
classical  Hardy-Littlewood inequality. \par

\vspace{3mm}

{\bf Remark 6.1.} Evidently, the {\it lower  bound }  for the "constant" $ W_{A,B}(p,q; r,d) $ in (6.4) is simultaneously lower bound
for arbitrary trace monomial Sobolev's inequality, not only for the radial functions. \par

\vspace{3mm}

\section{ Continuity of a function from weighted Sobolev-Grand Lebesgue spaces.}

\vspace{3mm}

We suppose in this section  $  p > D := D(A). $\par

  The following inequality was done again in the article  \cite{Cabre1} (Theorem 1.6):
there exists a constant $ C = C(D(A),p), $ depending only on $ p $ and $  D(A) $ such that
up to redefinition of the function $ u = u(x) $  in the set of zero measure

$$
\sup_{x \ne y} \left[ \frac{|u(x) - u(y)|}{|x-y|^{1-D(A)/p}}  \right] \le C(D(A),p) \cdot
\left( \int |x|^A \ |\nabla u(x)|^p  \ dx  \right)^{1/p}, \eqno(7.1)
$$
weighted version of the Morrey inequality, which may be rewritten as follows:

$$
\omega(u, \delta) \le C(D,p) \ \delta^{1 - D/p} \ \left( \int |x|^A \ |\nabla u(x)|^p  \ dx  \right)^{1/p}, \eqno(7.1a)
$$
where $ \omega(u, \delta), \ \delta > 0  $ is module of continuity  of the function $  u(\cdot). $ \par
 The  value $ C(D,p) = C(D(A),p) $ may be estimated after some calculation as follows:

$$
C(D,p) \le C_2(D) \cdot \frac{p}{p - D}, \ C_2(D) < \infty; \eqno(7.2)
$$
therefore

$$
\omega(u, \delta) \le C_2(D) \   \frac{p}{p - D} \cdot
\delta^{1 - D/p} \cdot \left( \int |x|^A \ |\nabla u(x)|^p  \ dx  \right)^{1/p}. \eqno(7.3)
$$

\vspace{3mm}

 Suppose now that the function $ |\nabla u| $ belongs to some GLS $  G(\psi) $ relative the measure $ \mu_A $ in the whole space $ R^m $
such that $  \supp \psi \subset (D,\infty): $

$$
\left( \int |x|^A \ |\nabla u(x)|^p  \ dx  \right)^{1/p} \le \psi(p) \ ||\nabla u||G(\psi_A), \ p \in \supp \psi.\eqno(7.4)
$$

 Denote

$$
\psi^{(D)}(p) = C_2(D) \cdot \frac{p}{p-D} \cdot \psi(p), \eqno(7.5)
$$
then it follows   from (7.3)

$$
\frac{\omega(u,\delta)}{ ||\nabla u||G(\psi^{(D)}_A) \cdot \delta }  \le \frac{\delta^{-D/p}}{1/\psi^{(D)}(p)}. \eqno(7.6)
$$
 Since the last inequality is true for all the values $  p  $ from the set $  \supp \psi = \supp \psi^{(D)},  $ we can take the minimum on
the right - hand side  of inequality (7.6):

$$
\frac{\omega(u,\delta)}{ ||\nabla u||G(\psi^{(D)}_A) \cdot \delta }  \le \inf_{p \in \supp \psi} \left[\frac{\delta^{-D/p}}{1/\psi^{(D)}(p)} \right] =
$$

$$
\left\{ \frac{1}{ \sup_{p \in \supp \psi}  [ \delta^{D/p} /\psi^{(D)}(p) ]}  \right\} =
\frac{1}{\phi( G(\psi^{(D)}), \delta^D  )}; \eqno(7.7)
$$
recall that $ \phi(G(\psi), \delta) $ denotes the fundamental function of the space $ G(\psi) $ at the point $ \delta. $ \par
 Thus, we proved in fact the following proposition. \par

 \vspace{3mm}

{\bf Theorem 7.1.} {\it  We obtained under formulate above conditions  that
up to redefinition of the function $ u = u(x) $  in the set of zero measure}

$$
\omega(u,\delta) \le \frac{||\nabla u||G(\psi^{(D)}_A) \cdot \delta }{ \phi( G(\psi^{(D)}), \delta^D  ) }, \ \delta > 0. \eqno(7.8)
$$

\vspace{3mm}

\section{ Concluding remarks.}

\vspace{3mm}

{\bf A.  Open question.} It is very interest by our opinion to find the sharp value of "constants"
for monomial weighted Sobolev's trace  inequality. \par

\vspace{3mm}

{\bf B. Second open question.} Is the Sobolev monomial inequality true for the {\it different powers} $ A,B $ and when $  p,q $
satisfying the relation (2.1) of theorem 2.1?

\vspace{3mm}

{\bf C. Extension.} It is no hard to generalise the monomial trace Sobolev's inequality to Grand Lebesgue Spaces,
but without exact constant computation. \par
 More interest problem is the extension of these inequalities on Lorentz, Marcinkiewicz etc. spaces, instead the
 classical Lebesgue-Riesz spaces.\par


\vspace{3mm}

\begin{thebibliography}{99}

\vspace{3mm}

\bibitem{Beesack1}
{\sc Beesack P.R.} Hardy's inequality and its extension. Pacific J. Math.;
{\bf 11}, (1961), 39-61.
\bibitem{Bennet1}
 {\sc Bennet G., Sharpley R.} Interpolation of operators. Orlando, Academic Press Inc., (1988).
\bibitem{Besov1}
{\sc O.V.Besov, P.I.Il'in, S.M. Nikol'skii.} Integral Representations of Functions
and Imbedding Theorems. Volume 1, (1978), John Wiley and Sons, Washington D.C., New York, Toronto, London, Sydney; A Halisted Press Book.
\bibitem{Besov2}
{\sc O.V.Besov, P.I.Il'in, S.M. Nikol'skii.} Integral Representations of Functions
and Imbedding Theorems. Volume 2, (1979), John Wiley and Sons, Washington D.C., New York, Toronto, London, Sydney; A Halisted Press Book.
\bibitem{Biezuner1}
{\sc R.J.Biezuner.} Best constants in Sobolev trace inequalities. Nonlinear Analysis,
{\bf 54}, (2003), 457-502.
\bibitem{Bradley1}
{\sc J.S.Bradley.}  Hardy inequalities with mixed norms. Canadian Math. Bull., 21(1978), p. 405-408.
\bibitem{Druet1}
{\sc O.Druet.} The best constants problem in Sobolev trace inequalities. Math. Annalen,
{\bf 314}, (1999), 327-346.
\bibitem{Edmunds1}
{\sc D.E.Edmunds and W.D. Evans.}  Sobolev Embeddings and Hardy Operators. In:
 Vladimir Maz'ya (Editor), "Sobolev Spaces in Mathematics", Part 1, International
 Mathematical Series,  Volume 8, Springer Verlag, Tamara Rozhkovskaya Publisher;
 (2009), New York, London, Berlin; p.153-184.
\bibitem{Escobar1}
{\sc J.F.Escobar.} Sharp constant in a Sobolev trace inequality. Indiana Math, J.,
{\bf 37}, (1988), 687-698.
\bibitem{Fiorenza1}
   {\sc Capone C., Fiorenza A., Krbec M.} On the Extrapolation Blowups in the
   $ L_p $ Scale. Collectanea Mathematica, {\bf 48}, 2, (1998), 71-88.
 \bibitem{Fiorenza2}
 {\sc A.Fiorenza.} Duality and reflexivity in grand Lebesgue spaces.
       Collectanea Mathematica (electronic version), {\bf 51}, 2, (2000), 131-148.
\bibitem{Fiorenza3}
 {\sc A. Fiorenza and G.E. Karadzhov.} Grand and small Lebesgue spaces and
       their analogs. Consiglio Nationale Delle Ricerche, Instituto per le
      Applicazioni del Calcoto Mauro Picine", Sezione di Napoli, Rapporto tecnico n.
      272/03, (2005).
\bibitem{Iwaniec1}
   {\sc T.Iwaniec and C. Sbordone.} On the integrability of the Jacobian under
      minimal hypotheses. Arch. Rat.Mech. Anal., 119, (1992), 129–143.
\bibitem{Iwaniec2}
 {\sc T.Iwaniec, P. Koskela and J. Onninen.} Mapping of finite distortion:
   Monotonicity and Continuity.  Invent. Math. 144 (2001), 507-531.
\bibitem{Kantorovicz1}
{\sc L.V.Kantorovicz, G.P.Akilov.} Functional Analysis. (1987) Kluvner Verlag.
\bibitem{Kozachenko1}
 {\sc Kozachenko Yu. V., Ostrovsky E.I.} (1985). The Banach Spaces of
      random Variables of subgaussian type. {\it Theory of Probab. and Math.
      Stat.} (in Russian). Kiev, KSU, {\bf 32}, 43-57.
\bibitem{Talagrand1}
 {\sc Ledoux M., Talagrand M.} (1991) Probability in Banach Spaces.
  Springer, Berlin, MR 1102015.
\bibitem{Lieb1}
{\sc E.H.Lieb.} Sharp constants in the Hardy-Littlewood-Sobolev and related inequalities.
Ann. Math., {\it 118},  (1983), 349-374.
\bibitem{Maz'ja1}
{\sc V.Maz'ja.} Sobolev Spaces. Kluvner Academic Verlag, (2002),
Berlin-Heidelberg-New York.
\bibitem{Mitrinovich1}
{\sc D.S.Mitrinovich, J.E. Pecaric and A.M.Fink.} Inequalities Involving
Functions and Their Integrals and Derivatives. Kluvner Academic Verlag,
(1996), Dorderecht, Boston, London.
\bibitem{Nazaret1}
{\sc B.Nazaret.} Best constant in Sobolev trace inequalities on the half space.
arXiv:math.FA/2109077 v1 21Aug 2005.
\bibitem{Ostrovsky1}
{\sc E.I. Ostrovsky.}  Exponential Estimations for Random Fields.
Moscow - Obninsk, OINPE, 1999 (Russian).
\bibitem{Ostrovsky2}
{\sc E. Ostrovsky and L.Sirota.} Moment Banach spaces: theory and applications.
HAIT Journal of Science and Engineering, {\bf C}, Volume 4, Issues 1-2,
pp. 233 - 262, (2007).
\bibitem{Ostrovsky3}
{\sc E.Ostrovsky, E.Rogover and L.Sirota.} Riesz's and Bessel's operators in in
bilateral Grand Lebesgue Spaces.
arXiv:0907.3321 [math.FA] 19 Jul 2009.
\bibitem{Ostrovsky4}
{\sc E. Ostrovsky and L.Sirota.} Weight Hardy-Littlewood inequalities for different
powers. arXiv:09010.4609v1[math.FA] 29 Oct 2009.
\bibitem{Ostrovsky5}
{\sc E. Ostrovsky E.} Bide-side exponential and moment inequalities
     for tail of distribution of Polynomial Martingales. Electronic
     publication, arXiv: math.PR/0406532 v.1 Jun. 2004.
\bibitem{Ostrovsky6}
{\sc E.Ostrovsky,  E.Rogover and L.Sirota.} Integral Operators in Bilateral Grand
Lebesgue Spaces. arXiv:09012.7601v1 [math.FA] 16 Dez. 2009.
\bibitem{Ostrovsky7}
{\sc E.Ostrovsky, L.Sirota.} Nikolskii-type inequalities for rearrangement invariant
spaces. arXiv:0804.2311v1 [math.FA] 15 Apr 2008.
\bibitem{Ostrovsky8}
 {\sc Ostrovsky E., Sirota L.} Moment Banach Spaces: Theory and Applications.
 HIAT Journal of Science and Engineering, Holon, Israel, v. 4, Issue 1-2, (2007), 233-262.
\bibitem{Sobolev1}
{\sc Sobolev S.L.} Some Applications of Functional Analysis into Mathematical Physic,
(1950), Publishing House LSU (Leningrad State University), (in Russian).

\bibitem{Stein1}
{\sc E.M.Stein.} {\it Singular Integrals and Differentiability Properties of Functions.} Princeton University
Press, Princeton, (1992).


\bibitem{Talenti1}
{\sc G.Talenti} Inequalities in Rearrangement Invariant Function Spaces. Nonlinear Analysis, Function
 Spaces and Applications. Prometheus, Prague, {\bf 5}, (1995), 177-230.
\bibitem{Young1}
{\sc Young Ja Park.} Sobolev Trace Inequalities.
arXiv:math.CA/0107065 v1 9Jul 2001.
\bibitem{Zhu1}.
{\sc M.Zhu} Some general forms of sharp Sobolev inequalities. J. Func. Anal., {\bf 156}, (1998), 75-120.

\vspace{11mm}


\bibitem{Bliss1}
{\sc G.A. Bliss.} {\it An integral inequality. }  J. London Math. Soc., 5 (1930), 40–46.

\bibitem{Cabre1}
{\sc Xavier Cabre and Xavier Ros-Oton. } {\it Sobolev and Isoperimetric Inequalities with Monomial Weight.}
arXiv:1210.4487v1 [math.AP] 16 Oct 2012

\bibitem{Ivanov1}
{\sc S. Ivanov, A. Nazarov. }  {\it On weighted Sobolev embedding theorems for functions with symmetries.}
St. Petersburg Math. J. 18 (2007), 77-88.

\bibitem{Ostrovsky500}
 {\sc Ostrovsky E., Sirota L.}
{\it Module of continuity for the functions belonging to the Sobolev-Grand Lebesgue Spaces.}
arXiv:1006.4177v1 [math.FA] 21 Jun 2010

\bibitem{Ostrovsky501}
 {\sc Ostrovsky E., Sirota L.} {\it Exact constant in Sobolev's and Sobolev's trace inequalities
for Grand Lebesgue Spaces.}
arXiv:1002.4865v1 [math.FA] 25 Feb 2010

\bibitem{Ostrovsky502}
 {\sc Ostrovsky E., Sirota L.} {\it Boundedness of operators in bilateral grand Lebesgue spaces,
with exact and weakly exact constant calculation.}
arXiv:1104.2963v1 [math.FA] 15 Apr 2011


\bibitem{Ostrovsky503}
 {\sc Ostrovsky E., Sirota L.} {\it Continuity of functions belonging to the fractional order Sobolev-Grand Lebesgue spaces.}
arXiv:1301.0132v1 [math.FA] 1 Jan 2013




\end{thebibliography}
\end{document}